\documentclass{article}[16pt]
\usepackage{amsfonts,amsmath}
\usepackage[russian]{babel}
\usepackage[cp1251]{inputenc}

\begin{document}

\title{Приближение линейных функционалов\\ на пространстве с выпуклой мерой}

\author{Д.В. Фуфаев \thanks{Механико-математический факультет Московского Государственного Университета имени М.В. Ломоносова,
Россия; e-mail:
fufaevdv@rambler.ru}  }

\date{}

\maketitle

Пусть $(X, \mu)$ - топологическое векторное пространство с заданной на нем вероятностной мерой. Существует два опеределения измеримого линейного функционала на $X$:

1) $f$ как отображение является измеримой и линейной функцией со значениями в $\mathbb R$ (иногда такие функционалы называют почти линейными).

2) $f$ - предел почти всюду непрерывных линейных функционалов.


В книге Гихмана и Скорохода [3, теорема VIII.1.1] было сформулировано, что для гильбертовых пространств эти два определения эквивалентны всегда, однако вскоре были построены контрпримеры, показывающие, что это неверно (например, в [4]). Поэтому актуальным стал вопрос - для каких классов мер на каких пространствах эти определения эквивалентны, то есть измеримый и линейный функционал приближается почти всюду конечномерными? Известно, что эти определения эквивалентны для класса центрированных гауссовских мер(см. [1, теорема 2.7.8]).

Есть гипотеза, что это верно для класса выпуклых мер. Класс выпуклых мер является обобщением класса гауссовских мер. Пока что выясняется, что свойства гауссовских мер присущи и выпуклым мерам, однако это выяснено далеко не для всех свойств, и соответствующие доказательства на порядок сложнее.

В работе исследуются  линейные функционалы для выпуклых продакт-мер на счетном произведении конечномерных пространств. Устанавливается эквивалентность этих двух определений для некоторых классов таких мер. Основным инструментом доказательства являются условные математические ожидания.

\begin{center}{\bf 1. Предварительные сведения.}\end{center}

\textbf{Определение~1.} {\it  Вероятностная мера $\mu$ на $\sigma$-алгебре $\sigma(X)$ в локально выпуклом пространстве $X$ называется выпуклой (или логарифмически вогнутой), если для всех множеств $A,B\in\sigma(X)$ и всех $\lambda\in[0,1]$ выполнено неравенство}

$$
\mu(\lambda A+(1-\lambda)B)\ge \mu(A)^{\lambda}\mu(B)^{1-\lambda}.
$$

Так как рассматриваемые меры будут произведениями конечномерных, то нам понадобится следующая характеризация конечномерных выпуклых мер:

\textbf{Лемма~1.} {\it  Мера $\mu$ на конечномерном пространстве $\mathbb R^m$ выпукла тогда и только тогда, когда найдутся аффинное отображение $T: \mathbb R^l\to \mathbb R^n$ и вероятностная мера $\nu$ с плотностью $p=e^{-V}$ на $\mathbb R^l$, где $V$ --- выпуклая функция (вне области определения $V$ полагаем $p=0$), причем $\nu=\mu \circ  T^{-1}$.
}

Отсюда, в частности, следует, что носитель меры $\mu$ в $\mathbb R^m$ - выпуклое множество, лежащее в некотором аффинном подпространстве.

Также нам понадобится следующее свойство выпуклых мер:

\textbf{Лемма~2.} {\it Пусть $X$ --- локально выпуклое пространство с выпуклой мерой $\mu$, $p$ --- измеримая относительно $\mu$ полунорма, конечная почти всюду. Тогда $p\in L^1(\mu).$
}

Учитывая, что если $f$ --- почти линейный измеримый функционал на $X$, то $|f|$ --- почти всюду конечная измеримая полунорма, получаем, что $f\in L^1(\mu).$

Доказательства лемм 1-2 можно найти в [5].

Напомним теперь определение условного математического ожидания. Пусть $(X, \mathfrak F, \mu)$ --- пространство с мерой и $\mathfrak B$ --- некоторая под-$\sigma$-алгебра в $\mathfrak F$.

\textbf{Определение~2.} {\it Пусть $f\in L^1 (\mu)$. Условным математическим ожиданием $f$ относительно $\sigma$-алгебры $\mathfrak B$ и меры $\mu$ называется такая $\mathfrak B$-измеримая $\mu$-интегрируемая функция $\mathbb E^{\mathfrak B}_{\mu}f$, что

$$
\int\limits_{X}gfd\mu=\int\limits_{X}g\mathbb E^{\mathfrak B}_{\mu}fd\mu
$$
для всякой ограниченной $\mathfrak B$-измеримой функции $g$.

}

Рассмотрим теперь пространство  $ \mathbb{R}^{\infty}$ как произведение конечномерных пространств:  $ \mathbb{R}^{\infty} =\prod\limits_{k=1}^{\infty}\mathbb{R}^{m_{k}}  $, где $\{m _{k}  \}_{k=1}^{\infty}$ - последовательность натуральных чисел. Пусть $\{\mu_{k_{}}\}_{k=1}^{\infty}$ - последовательность выпуклых мер, заданных на соответствующих $\mathbb{R}^{m_{k}}$. Тогда мера $\mu:=\bigotimes\limits_{k=1}^{\infty}     \mu_{k_{}}$, заданная на  $ \mathbb{R}^{\infty}$, также выпукла. Будем рассматривать именно такие меры.

Также введем обозначения $\tilde m_n=\sum\limits_{k=1}^{n}{m_{k}},\ \tilde \mu_{n}=\bigotimes\limits_{k=1}^{n}\mu_{k}.$  

\textbf{Лемма~3.} {\itПусть $\mathfrak{B}_{n}$ --- сигма-алгебра, порожденная первыми   $\tilde m_n$ координатными функциями. Тогда

$$
\mathbb {E}^{\mathfrak{B}_{n}}f(x_{1},\dots,x_{m_{n}})=
\int_{\prod\limits_{k=n+1}^{\infty}\mathbb{R}^{m_{k}}} f(x_{1},\dots,x_{m_{n}},x_{m_{n}+1},\dots)  \left(\bigotimes_{k=n+1}^{\infty}\mu_{k}\right) (d(x_{m_{n}+1},\dots)),
$$
где интегрирование ведется по произведению прямых, \newlineсоответствующему переменным $x_k$ с $k\ge m_n+1.$}

При этом, когда $f$ --- линейный функционал, верно следующее представление:

$$
\mathbb {E}^{\mathfrak{B}_{n}}f(x_{1},\dots,x_{m_{n}})=f(x_{1},\dots,x_{m_{n}},0,\dots) +  
$$

$$
  +\int_{\prod\limits_{k=n+1}^{\infty}\mathbb{R}^{m_{k}}}  f(0,\dots,0,x_{m_{n}+1},\dots)  \left(\bigotimes_{k=n+1}^{\infty}\mu_{k}\right) (d(x_{m_{n}+1},x_{m_{n}+2},\dots)).
$$
Этот функционал --- конечномерный аффинный (т.е. сумма линейного и константы), поэтому он непрерывен.

\textbf{Лемма~4.} {\it Функции $\mathbb {E}^{\mathfrak{B}_{n}}f$, рассматриваемые как функции на $\mathbb R^\infty$, сходятся к $f$ почти всюду при $n\to\infty$.}

Доказательства лемм 3-4 можно найти в [2, \S 10.1 и 10.2].

\begin{center}{\bf 2. Основные результаты.}\end{center}

\textbf {Теорема 1.} Пусть  найдется такой номер $q$, что носитель меры $\mu_q$ лежит в аффинном подпространстве, то есть в не проходящем через ноль. Тогда измеримый линейный функционал $f$ приближается непрерывными линейными $\mu$-п.в.

\textbf {Доказательство.}

Носитель меры $\mu_k$ --- это аффинное подпространство $H_k=  L_k+h_k$, где $L_k$ - линейное подпространство в $\mathbb{R}^{m_{k}}$ и $h_k \in H_k$.
Напомним, что найдется такой номер $q$, что носитель меры $\mu_q$ является аффинным, то есть не проходит через ноль (что равносильно $h_q\ne 0$). Тогда и конечномерные проекции $\mu$ на пространства $\mathbb{R}^{\tilde m_k}$  (т. е. меры $\tilde \mu_{k}=\bigotimes\limits_{j=1}^{k}\mu_{j}$ )  при $k\ge q$ также будут иметь аффинные носители вида $\tilde H_k=  \tilde L_k+\tilde h_k$, где $\tilde L_k=L_1\oplus L_2\oplus\dots\oplus L_k$, $\tilde h_k=h_1+ h_2+\dots+ h_k$, причем $\tilde h_k\ne 0$. 

Для $k \ge q$ построим линейные функционалы $\psi_k(x)$ на $\mathbb{R}^{\infty}$, такие что 

$$
\psi_k(\tilde h_k)=    \int_{\prod_{j=k+1}^{\infty}\mathbb{R}^{m_{j}}} f(0,\dots,0,x_{m_{k}+1},\dots)  \left(\bigotimes_{j=k+1}^{\infty}\mu_{j}\right) (d(x_{m_{k}+1},\dots))       ,
$$

$$
\psi_k|_{\tilde L_k}=0,\ \psi_k|_{\prod_{j=k+1}^{\infty}\mathbb{R}^{m_{j}}}=0
$$ 

(такие существуют, например, по теореме Хана-Банаха). После этого построим функционалы (уже аффинные) $\phi_k$ по формуле $\phi_k(x)=\psi_k(x-\tilde h_k)$, которые тоже будут зависеть только от первых $\tilde m_k$ координат. Теперь рассмотрим последовательность функционалов $g_k=\mathbb {E}^{\mathfrak{B}_{k}}f+\phi_k $. Это, очевидно, аффинные функционалы. Нетрудно проверить, что $g_k$ --- линейные функционалы (т.к. в начале координат принимают значение ноль). Кроме того, $g_k$ совпадает с $\mathbb {E}^{\mathfrak{B}_{k}}f$ на носителе меры  $\bigotimes\limits_{j=1}^{k}\mu_{j}$ (т.е. на $\tilde H_k$), поэтому они совпадают $\bigotimes\limits_{j=1}^{k}\mu_{j}$-почти всюду, а значит, так как оба зависят только от первых $\tilde m_k$ координат, $\mu$-почти всюду. Поэтому, так как $\mathbb {E}^{\mathfrak{B}_{k}}f$ сходились к $f$ $\mu$-почти всюду, то и $g_k$ будут сходиться к $f$ $\mu$-почти всюду. Отсюда заключаем, что $g_k$ --- последовательность конечномерных линейных функционалов, сходящихся к $f$ почти всюду.

\hfill$\Box$

\textbf {Теорема 2.} Пусть  носитель каждой меры $\mu_k$ симметричен относительно нуля (частный случай --- когда мера симметрична относительно начала координат). Тогда измеримый линейный функционал приближается конечномерными линейными $\mu$-п.в.


\textbf {Доказательство.} 

Пусть $\mu^-$ --- мера, полученная из $\mu$ отражением относительно нуля, она также будет произведением отражений исходных конечномерных мер, то есть $\mu^-=\bigotimes\limits_{k=1}^{\infty}     \mu^{-}_{k_{}}$, и все отражения также будут выпуклыми. Исходные меры $\mu, \mu_k$ будем обозначать как $\mu^+, \mu^+_k$. Функционал $f$ будет также $\mu^-$-измерим, поэтому можно построить его условные математические ожидания и по этой мере. Условные математические ожидания относительно $\mu^+, \mu^-$ будем обозначать $\mathbb {E}_+^{\mathfrak{B}_{n}}f, \mathbb {E}_-^{\mathfrak{B}_{n}}f$ соответственно.

Эти последовательности сходятся к $f$ почти всюду по мерам $\mu^+$ и $\mu^-$ соответственно, каждая по своей, вообще говоря. Но по условию, носитель каждой конечномерной меры симметричен относительно нуля, поэтому носитель отраженной меры будет совпадать с носителем исходной меры. Значит, исходная мера и отраженная будут абсолютно непрерывны относительно друг друга, значит, эти последовательности сходятся к $f$ почти всюду по обоим мерам.
Возьмем полусумму этих последовательностей, которая также будет сходиться к $f$ $\mu^+$-п.в. Очевидно, полусумма имеет следующий вид:

$$
f(x_{1},\dots,x_{m_{n}},0,\dots) +
 \frac{1}{2}  \int_{\prod\limits_{k=n+1}^{\infty}\mathbb{R}^{m_{k}}} f(0,\dots,x_{m_{n}+1},\dots)  \left(\bigotimes_{k=n+1}^{\infty}(\mu^+_{k}  +   \mu^-_{k})\right) (d(x_{m_{n}+1},\dots))
$$

Во втором слагемом в интеграле мера уже симметрична. Поэтому, так как интегрируется антисимметричная функция по симметричной мере, второе слагаемое равно нулю. Получили, что данная полусумма --- последовательность конечномерных линейных функций, сходящихся к $f$ $\mu$-п.в.

Отметим, что в случае симметричной меры отраженная совпадает с исходной, и если все $\mu_k$ были симметричными мерами, то и мера $\mu$ была бы симметричной, откуда следует, что искомый вид имели бы уже сами математические ожидания, и необходимости брать полусумму не имелось бы.

\hfill$\Box$

Из определения интеграла Лебега как верхней грани интегралов простых функций видно, что 

$$
\int_{\prod\limits_{j=k+1}^{\infty}\mathbb{R}^{m_{j}}} f(0,\dots,0,x_{m_{k}+1},\dots)  \left(\bigotimes\limits_{j=k+1}^{\infty}\mu^{+}_{j}\right) (d(x_{m_{k}+1},x_{m_{k}+2},\dots))=
$$

$$
=-\int_{\prod\limits_{j=k+1}^{\infty}\mathbb{R}^{m_{j}}} f(0,\dots,0,x_{m_{k}+1},\dots)  \left(\bigotimes\limits_{j=k+1}^{\infty}\mu^{-}_{j}\right) (d(x_{m_{k}+1},x_{m_{k}+2},\dots)).
$$

\textbf {Теорема 3.}
Если $f$ --- $\mu$-измеримый линейный функционал, причем числовая последовательность

$$
c_k= \int_{\prod\limits_{j=k+1}^{\infty}\mathbb{R}^{m_{j}}} f(0,\dots,0,x_{m_{k}+1},\dots)  \left(\bigotimes_{j=k+1}^{\infty}\mu_{j}\right) (d(x_{m_{k}+1},x_{m_{k}+2},\dots))
$$
стремится к нулю, то $f$ приближается непрерывными линейными функционалами $\mu$-п.в.

\textbf {Доказательство.} 

 Как и раньше, имеем

$$
\mathbb {E}_{\pm}^{\mathfrak{B}_{n}}f(x_{1},\dots,x_{m_{n}})=f(x_{1},\dots,x_{m_{n}},0,\dots) \pm c_n
$$

Тогда, для любого $x\in \mathbb{R}^{\infty} $ имеем

$$
|\mathbb {E}_-^{\mathfrak{B}_{n}}f(x) - f(x)| \le |\mathbb {E}_-^{\mathfrak{B}_{n}}f(x) - \mathbb {E}_+^{\mathfrak{B}_{n}}f(x)| + |\mathbb {E}_+^{\mathfrak{B}_{n}}f(x) - f(x)|
= 2 c_n + |\mathbb {E}_+^{\mathfrak{B}_{n}}f(x) - f(x)|
$$

Последнее выражение стремится к нулю для почти всех $x$ по мере $\mu$, так как  $\mathbb {E}_+^{\mathfrak{B}_{n}}f$ сходится к $f$ почти всюду, поэтому последовательность функционалов $\mathbb {E}_-^{\mathfrak{B}_{n}}f$ также сходится к $f$ почти всюду. Теперь, как и в предыдущем утверждении, взяв полусумму этих последовательностей, получаем искомую последовательность.

\hfill$\Box$

\begin{center}{\bf 3. Заключение.}\end{center}

Полученные результаты показывают эквивалентность определений измеримого линейного функционала для выпуклых продакт-мер на $\mathbb R^\infty$, кроме случая, когда носитель каждой меры лежит в линейном подпространстве $\mathbb R^{m_k}$, при этом нельзя найти аффинное подпространство, в котором он полностью содержится, и свойства навроде симметричности носителя нет. Теорема 3 дает достаточное условие для эквивалентност, поэтому она также дает один из возможных путей обработки оставшегося случая --- доказать, что условие теоремы 3 в этом случае будет выполняться.

\end{document}